\documentclass[12pt, reqno]{amsart}
\usepackage{amsmath, amstext, amsbsy, amssymb, amscd}
\usepackage{amsmath}
\usepackage{amsxtra}
\usepackage{amscd}
\usepackage{amsthm}
\usepackage{amsfonts}
\usepackage{amssymb}
\usepackage{eucal}
\usepackage{color}

\newtheorem{theorem}{Theorem}[section]
\newtheorem{lemma}{Lemma}[section]
\newtheorem{definition}{Definition}[section]
\newtheorem{prop}{Proposition}[section]

\theoremstyle{remark}
\newtheorem{remark}{Remark}[section]

\newcommand{\C}{\mathbb C}
\newcommand{\cinf}{ c_{\infty} }
\newcommand{\dinf}{ d_{\infty} }
\newcommand{\gl}{ {\mathfrak {gl}} }
\newcommand{\hgl}{ \widehat{ \mathfrak{gl} } }
\newcommand{\glpm}{ \widehat{{ \mathfrak {gl} }}_{\pm} }
\newcommand{\hf}{ \frac12}
\newcommand{\half}{ \hf}
\newcommand{\hL}{{\Lambda} }
\newcommand{\vac}{ |0 \rangle }

\newcommand{\la}{\lambda}
\newcommand{\Z}{\mathbb Z}
\newcommand{\HZ}{\frac{1}{2}+\Z}
\newcommand{\fock}{\mathcal F}
\newcommand{\trace}{\,{\rm tr}\,}

\newcommand{\PA} {\mathcal P }

\newcommand{\ep}{\varepsilon}
\newcommand{\eps}{\epsilon}
\newcommand{\no} {\text{:}}

\newcommand{\vphi}{\varphi}
\newcommand{\A}{\mathsf A}
\newcommand{\D}{\mathsf D}

\begin{document}

\title[Bloch-Okounkov functions: a Classical half-integral case]
{The Bloch-Okounkov correlation functions, a Classical half-integral case}

\author{David G. Taylor}
\address{Department of Mathematics, Computer Science, and Physics, Roanoke College, Salem, VA 24153}
\email{taylor@roanoke.edu}

\begin{abstract}
Bloch and Okounkov's correlation function on the infinite wedge space has connections to Gromov-Witten theory, Hilbert schemes, symmetric groups, and certain character functions of $\hgl_\infty$-modules of level one.  Recent works have calculated these character functions for higher levels for $\hgl_\infty$ and its Lie subalgebras of classical type.  Here we obtain these functions for the subalgebra of type $D$ of half-integral levels and as a byproduct, obtain $q$-dimension formulas for integral modules of type $D$ at half-integral level.
\end{abstract}

\maketitle

%\tableofcontents


\section{Introduction}

Bloch and Okounkov \cite{BO} introduced an $n$-point correlation function on the
infinite wedge space and found an elegant closed formula in terms of
theta functions.  From a representation theoretic viewpoint, the Bloch-Okounkov
$n$-point function can be also easily interpreted as correlation
functions on integrable modules over Lie algebra $\hgl_\infty$ of
level one (cf. \cite{Ok, Mil, CW}). Along this line, Cheng and Wang
\cite{CW} formulated and calculated such $n$-point
correlation functions on integrable $\hgl_\infty$-modules of level
$l$ ($l \in \mathbb N$).  These correlation functions proved to be very useful in many applications such as in Gromov-Witten theory, Hilbert schemes, and the study of the symmetric groups.

The author and Wang \cite{TW} extended the formulation and computation of these correlation functions to the other classical subalgebras of $\hgl_\infty$; there we have calculated the $n$-point correlation functions for integrable modules of arbitrary positive level for the subalgebras classically identified as $b_\infty$, $c_\infty$, and $d_\infty$.  The author, along with Cheng and Wang \cite{CTW}, later further extended results to modules of negative level for $\hgl_\infty$ and its same subalgebras.  For more history of this problem, we refer the reader to the introduction of \cite{TW}.

It has been known since Hermann Weyl that representations of the orthogonal Lie algebras have certain annoying features due to two connected components of the orthogonal group; in particular, Weyl's character formula gives a better result for the character of certain pairs of representations rather than the individual components.  In \cite{TW}, we are forced to consider the direct sum of two irreducible $d_\infty$-modules for this reason.  In this paper, we aim to examine the case for $d_\infty$-modules of positive half-integral level.  Our main strategy, as in \cite{TW}, is to use a free-field realization \cite{DJKM2} and a Howe duality due to Wang \cite{W1} between $\dinf$ and the Lie group $O(2l+1)$.  We develop an operator in $d_\infty$ that is able to distinguish between the two components of this direct sum and use this operator to help compute a formula for the $n$-point correlation functions on the irreducible $\dinf$-modules.

The paper is organized as follows.  In section~\ref{sec:two}, we review some of the preliminaries.  First we review the definitions and notations we will use regarding $\hgl_\infty$ and $\dinf$.  Then we give a brief review of the Lie group $O(2l+1)$ and conclude with a quick review of the problem in the $\hgl_\infty$ case.  This section also introduces some of the Fock space definitions as well as the original Bloch-Okounkov function which will appear in several of our formulas.  Finally, in section~\ref{sec:four} we present our main theorems with proof.

\textbf{Acknowledgement} The author is partially supported by a faculty research grant from Roanoke College.  The author wishes to thank Weiqiang Wang for many helpful discussions and comments.  The author also wishes to thank the referee for helpful suggestions and prompt correspondence.

\section{The preliminaries}\label{sec:two}

\subsection{Classical Lie algebras of infinite dimension}
\label{sec_algebras}
In this subsection we review Lie algebras $\hgl = \hgl_\infty$
and Lie subalgebras of type $D$ (cf. \cite{DJKM2, K}).
\subsubsection{Lie algebra $\hgl$}

Denote by $\gl$ the Lie algebra of all matrices $(a_{ij})_{i,j \in
\Z}$ satisfying $a_{ij} =0$ for $|i -j|$ sufficiently large.
Denote by $E_{ij}$ the infinite matrix with $1$ at $(i, j)$ place
and $0$ elsewhere and let the weight of $ E_{ij}$ be $j - i $.
This defines a $\Z$--principal gradation $\gl = \bigoplus_{j \in
\Z} \gl_j$. Denote by $\hgl \equiv \hgl_\infty = \gl \oplus \C C$
the central extension given by the following $2$--cocycle with
values in $\C$:
\begin{eqnarray}
 C(A, B) = \trace \left([J, A]B \right)
  \label{eq_cocy}
\end{eqnarray}
where $J = \sum_{ j \leq 0} E_{ii}$. The $\Z$--gradation of Lie
algebra $\gl$ extends to $\hgl$ by letting the weight of $C$ to be
$0$. This leads to a triangular decomposition
$$\hgl = \widehat{\gl}_{+}
          \oplus \widehat{\gl}_{0} \oplus \widehat{\gl}_{-}  $$
where
$ \glpm = \oplus_{ j \in \mathbb N} \widehat{\gl}_{\pm j},
\widehat{\gl}_{0} = \gl_0 \oplus \C C.$
Let
\begin{eqnarray*}
 H^a_i  = E_{ii} - E_{i+1, i+1} + \delta_{i,0} C \quad (i \in \Z).
\end{eqnarray*}
Denote by $L(\hgl; \Lambda)$ the highest weight $\hgl$--module
with highest weight $\Lambda  \in {\hgl}_{0}^{*}$, where $C$ acts
as a scalar which is called the {\em level}. Let $\hL_j^a \in
\hgl_0^{*}$ be the fundamental weights, i.e. $\hL_j^a ( H_i^a ) =
\delta_{ij}.$

\subsubsection{Lie algebra $\dinf$}
  \label{subsec_dinf}
Let
\begin{eqnarray*}
{\overline{d}}_{\infty} = \{ (a_{ij})_{i,j \in \Z} \in \gl \mid
a_{ij} = -a_{1-j,1-i} \}
 \end{eqnarray*}
be a Lie subalgebra of $\gl$ of type $D$. Denote by $\dinf =
{\overline{d}}_{\infty} \bigoplus \C C $ the central extension
given by the $2$-cocycle (\ref{eq_cocy}). Then $\dinf$ has a
natural triangular decomposition induced from $\hgl$ with Cartan
subalgebra ${\dinf}_0 =  \widehat{\gl}_{0}  \cap \dinf$.  Given
$\Lambda \in {\dinf}_0^* $, we let
\begin{eqnarray*}
  H^d_i & = & E_{ii} + E_{-i, -i} - E_{i+1, i+1} - E_{-i+1, -i+1}
                                \quad (i \in \mathbb N),\\
  H^d_0 & = & E_{0,0} + E_{-1,-1} -E_{2,2} -E_{1,1} + 2C.
\end{eqnarray*}
Denote by $\hL^d_i $ the $i$-th fundamental weight of $\dinf$,
i.e. $\hL^d_i (H^d_j ) = \delta_{ij}$.

\subsection{Classical Lie group $O(2l+1)$}
    \label{sec_classical}

\subsubsection{$O(2l+1)$}

Let $O (2l+1) = \{ g \in GL(2l+1) \mid {}^t g J g = J \}$, where
  \begin{eqnarray*}
  J =\left[ \begin{array}{ccc}
    0      & I_l   & 0 \\
    I_l    & 0     & 0   \\
    0      & 0     & 1
         \end{array} \right].
  \end{eqnarray*}
The Lie algebra ${\mathfrak {so}}(2l+1)$ is the Lie subalgebra of
${\mathfrak {gl}}(2l +1)$ consisting of $(2l+1)\times (2l+1)$
matrices of the form
  \begin{eqnarray}
  \left[ \begin{array}{ccc}
    \alpha      & \beta         & \delta \\
    \gamma      & -{}^t{\alpha} & h      \\
     -{}^t h    & -{}^t{\delta} & 0
         \end{array} \right]
   \label{eq_form}
  \end{eqnarray}
where $\alpha, \beta, \gamma$ are $l \times l$ matrices and $
\beta, \gamma$ skew-symmetric. The Borel subalgebra ${\mathfrak b}
( {\mathfrak {so} } (2l+1) ) $ consists of matrices
(\ref{eq_form}) by putting $\gamma $, $h$, $\delta$ to be $0$ and
$\alpha$ to be upper triangular. The Cartan subalgebra ${\mathfrak
h} ( {\mathfrak {so} } (2l+1) ) $ consists of diagonal matrices of
the form $ \text{diag} ( t_1, \ldots, t_l; -t_1 \ldots -t_l; 0 ),$
$ t_i \in \C$. An irreducible module of $SO(2l+1)$ is
parameterized by its highest weight $(m_1, \ldots, m_l ) \in
\PA^l$, where $\PA^l$ denotes the set of partitions with at most
$l$ non-zero parts.

It is well known that $ O (2l+1)$ is isomorphic to the direct
product $ SO (2l+1) \times \Z_2$ by sending the minus identity
matrix to $-1 \in \Z_2 = \{ \pm 1 \}.$ Denote by $\det$ the
non-trivial one-dimensional representation of $O(2l+1)$. An
representation $\lambda$ of $SO(2l+1)$ extends to two different
representations $\lambda$ and $\lambda \otimes \det$ of $O(2l+1)$.
Then we can parameterize irreducible representations of $O(2l+1)$
by $ (m_1, \ldots, m_l)$ and $ (m_1, \ldots, m_l) \otimes \det $.
We shall denote
$$ \Sigma(B) =
 \PA^l \cup \left\{\la \otimes \det  \mid \la \in \PA^l \right \}.
 $$
For more details regarding a parametrization of irreducible
modules of various classical Lie groups including $O(2l+1)$, we refer the reader to \cite{BtD}.

\subsection{The Fock space $\fock^l$}

Consider a pair of fermionic fields
$$
 \psi^{+}(z) = \sum_{n \in \Z+\hf } \psi^{+}_n z^{ -n -\hf},
  \quad \psi^{-}(z) = \sum_{n \in \Z+\hf }
                           \psi^{-}_n z^{ -n -\hf},
 $$
with the following anti-commutation relations
\begin{eqnarray*}
 {[}\psi^{+} _m, \psi^{-}_n ]_{+}
%&=& \psi^{+}_m \psi^{-}_n + \psi^{-}_n \psi^{+}_m
= \delta_{m+n,0},  \qquad
 {[}\psi^{\pm} _m, \psi^{\pm}_n ]_{+} = {0}.
\end{eqnarray*}
Denote by $\mathcal F$ the Fock space of the fermionic fields
$\psi^{\pm}(z)$ generated by a vacuum vector $\vac$ which
satisfies
\begin{eqnarray*}
  \psi^{-}_n \vac = \psi^{+}_n \vac = 0
  & \text{ for } n \in \hf +\Z_+.
\end{eqnarray*}
We have the standard charge decomposition (cf. \cite{MJD})
$$ \fock =\bigoplus_{k \in \Z} \fock^{(k)}.
$$
Each $\fock^{(k)}$ becomes an irreducible module over a certain
Heisenberg Lie algebra. The shift operator $\textsf{S}:
\fock^{(k)} \rightarrow \fock^{(k+1)}$ matches the highest weight
vectors and commutes with the creation operators in the Heisenberg
algebra.

Now we take $l$ pairs of fermionic fields, $\psi^{\pm,p} (z) \;( p
= 1, \dots, l)$ and denote the corresponding Fock space by
$\fock^l$. Introduce the following generating series
\begin{eqnarray}
  E(z,w) & \equiv & \sum_{i,j \in \Z} E_{ij} z^{i }
     w^{-j  }
   = \sum_{p =1}^l \no\psi^{+,p} (z) \psi^{-,p} (w)\no,  \label{eq_genelin}
\end{eqnarray}
where the normal ordering $::$ means that the operators
annihilating $\vac$ are moved to the right with a sign. It is well
known that the operators $E_{ij} \; (i,j \in \Z)$ generate a
representation in $\fock^l$ of the Lie algebra $\hgl$ with level
$l$.

Let
\begin{eqnarray} \label{e-}
e^-_{pq}  =\sum_{r\in
\Z}\no\psi_{r}^{-p}\psi^{-q}_{-r}\no,\quad
e_{pq}^+ = \sum_{r\in\Z}
\no\psi_{r}^{+p}\psi^{+q}_{-r}\no, \quad p \neq q,
\end{eqnarray}
 and let
\begin{eqnarray}  \label{e+}
e_{pq}  =\sum_{r\in\Z}
\no\psi_{r}^{+p}\psi^{-q}_{-r}\no +\delta_{pq} \eps.
\end{eqnarray}
The operators $e^+_{pq}, e_{pq}, e^-_{pq} \;(p, q = 1, \cdots, l)$
generate Lie algebra $\mathfrak {so}(2l)$ (cf. \cite{FF, W1}).

\subsection{The main results of \cite{BO, CW}}

Recall that Bloch and Okounkov \cite{BO} introduced  the following
operators in $\hgl$
$$\no\A(t)\no
= \sum_{k\in\Z} t^{k-\hf}E_{k,k}, \quad
\A(t) = \no\A(t)\no + \frac{1}{t^\hf-t^{-\hf}} C.
$$
%We easily verify that
%
%\begin{equation} \label{eq:AD}
%\D(t) = \A(t) - \A(t^{-1}), \qquad \no \D(t)\no = \no
%\A(t)\no - \no \A(t^{-1})\no.
%\end{equation}

Given $\la =(m_1,\ldots,m_l)\in \Sigma(A)$, we denote by $\hL(\la)$
the $\hgl$-highest weight $\hL^a_{m_1}+\cdots+\hL^a_{m_l}$. The
energy operator $L_0$ on the $\hgl$-module $L(\hgl;\hL(\la))$ with
highest weight vector $v_{\hL(\la)}$ is characterized by
\begin{eqnarray}
L_0 \cdot v_{\hL (\la)}
 &=& \hf \Vert\la\Vert^2  \cdot v_{\hL (\la)}, \label{weight} \\
 {[}L_0, E_{ij}] &=& (i-j) E_{ij}, \nonumber
\end{eqnarray}
where $$\Vert\la\Vert^2 :=\la_1^2 +\la_2^2 +\cdots + \la_l^2,$$
On $\fock^l$, we can realize $L_0$ as
$$L_0 =
\sum_{p=1}^l\sum_{k\in\Z+\hf}
k\no\psi_{-k}^{+,p}\psi_{k}^{-,p}\no.$$

The {\em $n$-point $\hgl$-correlation function of level $l$}
associated to $\la$ is defined in \cite{BO} for $l=1$ and in
\cite{CW} for general $l$ as
\begin{equation*}
\mathfrak A_\la^l (q;{\bf t}) \equiv
 \mathfrak A_\la^l (q;t_1,\ldots,t_n)
 :=\trace_{L(\hgl,\hL(\la))}(q^{L_0} \A(t_1)\A(t_2)\cdots
 \A(t_n)).
\end{equation*}
Here and below we denote ${\bf t} =(t_1, \ldots, t_n)$.

Let $ (a; q)_\infty := \prod_{r=0}^\infty (1-aq^r)$. Define the
theta function
\begin{eqnarray}
 \label{theta:def} \Theta (t) &:=&
 (t^{\hf} -t^{-\hf})(q;q)_\infty^{-2} (qt; q)_\infty(qt^{-1};q)_\infty \\
 \label{theta:deriv} \Theta^{(k)} (t) &:=&
 \left(t\frac{d}{dt} \right)^k \Theta (t), \quad\text{ for }k\in\Z_+.
\end{eqnarray}
Denote by $F_{bo}(q; {\bf t})$ or $F_{bo}(q;t_1,\ldots,t_n)$ the
following expression
\begin{eqnarray}  \label{eq:bo}
\frac{1}{(q; q)_\infty}\cdot \sum_{\sigma\in S_n} \frac{{\rm det}
\Big{(}\frac{\Theta^{(j-i+1)}(t_{\sigma(1)}\cdots
t_{\sigma(n-j)})}{(j-i+1)!} \Big{)}_{i,j=1}^n}
{\Theta(t_{\sigma(1)}) \Theta(t_{\sigma(1)}t_{\sigma(2)})\cdots
\Theta(t_{\sigma(1)}t_{\sigma(2)}\cdots t_{\sigma(n)})}.
\end{eqnarray}
It is understood here that $1/ (-k)!=0$ for $k>0$, and for $n=1$,
we have $F_{bo}(q; t) = (q;q)_\infty^{-1} \Theta(t)^{-1}$. The
following summarizes the main results of Bloch-Okounkov \cite{BO}
for $l=1$ and Cheng-Wang \cite{CW} for general $l \ge 1$.
\begin{theorem}  \label{th:CWmain}
Associated to $\la =(\la_1,\ldots,\la_l)$, where $\la_1\ge\ldots
\ge \la_l$ and $\la_i \in\Z$, the $n$-point $\hgl$-function of
level $l$ is given by
\begin{eqnarray*}
\mathfrak A^l_\la(q; {\bf t}) =
 q^{\frac{\Vert\la\Vert^2}{2}} (t_1t_2\cdots t_n)^{|\la|} \prod_{1\le i<j\le
l} (1-q^{\la_i-\la_j+j-i}) \cdot
 F_{bo}(q;{\bf t})^l
\end{eqnarray*}
where $|\la| := \la_1 +\cdots +\la_l.$
\end{theorem}
In the simplest case, i.e. $l =n =1$, we have
$$
\mathfrak A^1_\la(q; t) =q^{\frac{\la^2}{2}}t^\la\cdot F_{bo}(q; t)
=\frac{q^{\frac{\la^2}{2}}t^\la}{(q;q)_\infty \Theta(t)}.
$$

\section{The correlation functions on $d_\infty$-modules}
\label{sec:four}

Let $t$ be an indeterminate and define the following operators in $d_\infty$:
\begin{eqnarray*}
 \no\D(t)\no &=& \sum_{k\in \mathbb N}(t^{k-\hf} -
t^{\hf-k})(E_{k,k}-E_{1-k,1-k}), \\
 \D(t) &=& \no \D(t)\no +\frac{2}{t^\hf-t^{-\hf}}C.
\end{eqnarray*}

\begin{definition}
The $n$-point $\dinf$-correlation function of level $l+\hf$
associated to $\la \in \PA^l \cup \PA^l\otimes\det$, denoted by $\mathfrak
D^{l+\hf}_\la(q,{\bf t})$ or also by $\mathfrak
D^{l+\hf}_\la(q,t_1,\ldots, t_n)$, is
\begin{eqnarray*}
\trace_{L(d_\infty;\Lambda(\la))} q^{L_0}\D(t_1)\cdots \D(t_n).
\end{eqnarray*}
\end{definition}

\begin{remark} In \cite{TW}, the trace is taken over the direct sum of the modules $L(\dinf;\Lambda(\la))$ and $L(\dinf;\Lambda(\la\otimes\det))$ for $\la$ in $\PA^l$ for technical reasons.\end{remark}

\subsection{Fock space $\fock^{l+\hf}$ and $\D(t)$ realization}

Consider the neutral fermion
$$
\vphi(z) = \sum_{n\in \Z+\hf}\vphi_nz^{-n-\hf}
$$
which satisfies the commutation relation
$$
[\vphi_m,\vphi_n]_+ = \delta_{m,-n}.
$$
We denote by $\fock^{l+\hf}$ the Fock space of one neutral fermion
$\vphi(z)$ and $l$ pairs of complex fermions $\psi^{\pm,p}(z),
1\le p \le l$, generated by a vacuum vector $\vac$ which satisfies
\begin{eqnarray*}
  \vphi_n \vac =\psi^{+,p}_n \vac = \psi^{-,p}_n \vac = 0
  & \text{ for } n \in \hf +\Z_+.
\end{eqnarray*}

Let
$$e^\pm_{p}  =\sum_{r\in \Z}\no \psi_{r}^{\pm p}\vphi_{-r}\no ,
 \quad 1\le p\le l.
$$
It is known (cf. \cite{FF, W1}) that the above operators $e^+_p,
e^-_p$ together with $e^+_{pq}, e_{pq}, e^-_{pq}$ $(p, q = 1,
\cdots, l)$ defined in (\ref{e-}, \ref{e+}) generate Lie algebra
$\mathfrak {so}(2l+1)$.

When acting on $\fock^{l+\half}$, we may then rewrite $\D(t)$ as $$\D(t) = \sum_{k\in\HZ}t^k \left(\sum_{i=1}^l
(\psi^{+,i}_{-k}\psi^{-,i}_{k} + \psi^{-,i}_{-k}\psi^{+,i}_k) +
\vphi_{-k}\vphi_{k}\right).$$

For later use, we have the following lemma giving an isomorphism of Fock spaces.

\begin{lemma} \label{lem:fockiso}
Given a pair of complex fermions $\psi^{\pm} (z)$, we let
$$
\varphi_n := (\psi^+_n +\psi^-_n)/\sqrt{2},
 \qquad \varphi_n' :=  i(\psi^+_n-\psi_n^-)/\sqrt{2}.
$$
Then, $\varphi_n$ and $\varphi_n'$ satisfy the anti-commutation
relations:
\begin{eqnarray*}
 {[}\varphi_n,\varphi_m]_+ = \delta_{n,-m}, &&
 {[}\varphi_n',\varphi_m']_+ = \delta_{n,-m},\\
 {[}\varphi_n,\varphi_m'{]}_+ = 0, && \text{for } m,n\in \underline{\Z}.
\end{eqnarray*}
Hence, there is an isomorphism of Fock spaces
\[\fock^\half\otimes\fock^\half \cong \fock\]
 \end{lemma}
\begin{proof}
Verified by a direct computation.
\end{proof}

\subsection{The $n$-point $\dinf$-correlation functions of level $l+\hf$}

Consider the $\dinf$ operator $$\alpha = \sum_{k>0} \varphi_{-k}\varphi_k$$ and set the following notation
$$\begin{aligned}
\mathbb
{D}^{l+\hf}_\la(q,{\bf t}) &= \trace_{L(d_\infty;\Lambda(\la)) \oplus L(d_\infty;\Lambda(\la
\otimes \det))} q^{L_0}\D(t_1)\cdots \D(t_n) \\
\overline{\mathbb
{D}}^{l+\hf}_\la(q,{\bf t}) &= \trace_{L(d_\infty;\Lambda(\la)) \oplus L(d_\infty;\Lambda(\la
\otimes \det))} \left(-1\right)^\alpha q^{L_0}\D(t_1)\cdots \D(t_n).
\end{aligned}$$

Note that $\mathbb{D}^{l+\hf}_\la(q,{\bf t})$ was computed in \cite{TW}.

\begin{prop}\cite[Theorem 4.1]{TW} \label{th:lplushalf}
The function $\mathbb D^{l+\hf}_\la(q,{\bf t})$ is equal to
\begin{eqnarray*}
 && \mathbb D^{\hf}_{(0)} (q;{\bf t}) \times \\
 &&\times
 \sum_{\sigma\in W(B_l)}\left(-1\right)^{\ell(\sigma)}
q^{\frac{\Vert \la+\rho-\sigma(\rho) \Vert^2}{2}}
\prod_{a=1}^l
 \Big(
\sum_{\vec{\eps}_a\in\{\pm 1\}^n}
 [\vec{\eps}_a] (\Pi {\bf t}^{\vec{\eps}_a})^{k_a}
 F_{bo}(q;{\bf t}^{{\vec{\eps}_a}})
 \Big)
\end{eqnarray*}
where $k_a = (\la+\rho-\sigma(\rho), \ep_a)$, $W(B_l)$ is the Weyl group of type $B$, $\rho$ is the usual half-sum of positive roots for type $B$, $\vec{\eps} =(\eps_1,\eps_2,\ldots,\eps_n)$,
$[\vec{\eps}] =\eps_1\eps_2 \cdots \eps_n$, and $\Pi{\bf
t}^{\vec{\eps}} =t_1^{\eps_1}\cdots t_n^{\eps_n}$.  Also, $\mathbb D^\hf_{(0)}(q;{\bf t})$ is given in \cite[Proposition 4.2]{TW}.
\end{prop}

The following formula from \cite{TW} will be used later.

\begin{equation}\label{needed}
 \trace_\fock z^{e_{11}}q^{L_0}\D(t_1)\cdots\D(t_n) = \sum_{k\in\Z} z^k q^{\frac{k^2}{2}}
 \sum_{\vec{\eps}  \in \{\pm 1\}^n}  [\vec{\eps}] \cdot
 \left(\Pi{\bf t}^{\vec{\eps}} \right)^k \cdot
 F_{bo}(q;{\bf t}^{\vec{\eps}})
\end{equation}

The main results of this paper are the computation of the function $\overline{\mathbb{D}}^{l+\half}_\la(q,{\bf t})$ and Proposition \ref{maintheorem} below.  A recursive formula for $\overline{\mathbb{D}}_{(0)}^\half(q,{\bf t})$ can be obtained similar to \cite[Proposition 4.2]{TW}.  Note that \cite[Theorem 8]{W2} implicitly gives the 1-point version as $$\overline{\mathbb{D}}_{(0)}^\half(q,t) = -(q^{\half})_\infty\left(\sum_{n=1}^\infty \frac{q^{n-\half}(t^{n-\half}-t^{-n+\half})}{1-q^{n-\half}}\right)+\frac{t^\half}{t-1}(q^\half)_\infty.$$

\begin{lemma} We have $$[\alpha,\varphi_r] = \varphi_r,\qquad\qquad [\alpha,\psi_r^\pm] = 0.$$ Equivalently, $\alpha$ acts on vectors of $d_\infty$-modules by counting the number of $\varphi_r$s in the vector.\end{lemma}

\begin{proof} \label{shortlemma} The lemma follows by direct computation using the anti-commutation relations amongst the $\varphi$s and $\psi$s.\end{proof}

\begin{prop}\label{maintheorem} For $\la\in\PA^l$, the $n$-point $\dinf$-correlation functions of level $l+\hf$ are given by $$\mathfrak{D}^{l+\half}_\la(q,{\bf t}) = \frac{\mathbb D^{l+\half}_\la(q,{\bf t}) + \overline{\mathbb{D}}^{l+\half}_\la(q,{\bf t})}{2}$$
$$\mathfrak{D}^{l+\half}_{\la\otimes\det}(q,{\bf t}) = \frac{\mathbb D^{l+\half}_\la(q,{\bf t}) - \overline{\mathbb{D}}^{l+\half}_\la(q,{\bf t})}{2}.$$
\end{prop}

\begin{proof}
Using Lemma \ref{shortlemma}, the operator $\alpha$ acting on an element of either $L(d_\infty;\Lambda(\la))$ or $L(d_\infty;\Lambda(\la
\otimes \det))$ counts the number of $\varphi$s in the vector.  The structure of the highest weight vectors for these modules is well-known (cf. \cite[Theorem 4.1]{W1}) and elements of $L(d_\infty;\Lambda(\la))$ (respectively $L(d_\infty;\Lambda(\la\otimes\det))$) have an even (respectively odd) number of $\varphi$s; hence $\left(-1\right)^\alpha$ acts as 1 on $L(d_\infty;\Lambda(\la))$ and as $-1$ on $L(d_\infty;\Lambda(\la\otimes \det))$ and the result follows.
\end{proof}

We set $$\alpha' = \sum_{k>0}\left(\varphi_{-k}\varphi_{k} + \varphi_{-k}'\varphi_k'\right)$$ and given a subset $I = (i_1,\ldots, i_s) \subseteq \{1,\ldots,n\}$, we
denote by $I^c$ the complementary set to $I$, and ${\bf t}_I
=(t_{i_1},\ldots,t_{i_s})$. By convention, we let
\begin{eqnarray} \label{eq:empty}
\overline{\mathbb D}^{\hf}_{(0)}(q,{\bf t}_\emptyset) =\trace_{\fock^{\hf}}
\left(-1\right)^{\alpha'} q^{L_0} =(q^\hf;q)_\infty.
\end{eqnarray}

\begin{prop} \label{recursive}
We have
\begin{eqnarray} \label{eq:subset}
%\bff (1,q;t_1,\ldots,t_n)
 \trace_\fock \left(-1\right)^{\alpha'} q^{L_0}\D(t_1)\cdots\D(t_n)
 &= \displaystyle\sum_{I\subseteq \{1,\ldots,n\}}
 \overline{\mathbb D}^{\hf}_{(0)}(q,{\bf t}_I)
 \overline{\mathbb D}^{\hf}_{(0)}(q,{\bf t}_{I^c}).
\end{eqnarray}

Equivalently, we have
$$
\begin{aligned}
\overline{\mathbb D}^{\hf}_{(0)} (q,{\bf t})
  &=& \hf (q^\hf;q)_\infty^{-1}
  \left( \sum_{k\in\Z}  \left(-1\right)^k q^{\frac{k^2}{2}}
 \sum_{\vec{\eps}  \in \{\pm 1\}^n}  [\vec{\eps}] \cdot
 \left(\Pi{\bf t}^{\vec{\eps}} \right)^k
 F_{bo}(q;{\bf t}^{\vec{\eps}}) \right. \\
 && -\left. \sum_{\emptyset \subsetneq I \subsetneq\{1,\ldots,n\}}
 \overline{\mathbb D}^{\hf}_{(0)}(q,{\bf t}_I)
 \overline{\mathbb D}^{\hf}_{(0)}(q,{\bf t}_{I^c}) \right).
\end{aligned}
$$
\end{prop}

\begin{proof}
A simple calculation reveals that
$$
\psi_{-k}^+\psi_k^-+\psi_{-k}^-\psi_k^+ = \varphi_{-k}\varphi_k +
\varphi_{-k}'\varphi_k'.
$$
so under the isomorphism $\fock \cong \fock^\half\otimes\fock^\half$, we may write $\D(t) = \D_1(t) + \D_2(t)$ where $\D_1(t) = \sum_{k\in\Z+\half}t^k \varphi_{-k}\varphi_k$ and
$\D_2(t) = \sum_{k\in\Z+\half} t^k\varphi_{-k}'\varphi_k'$.
Therefore, we have
$$\begin{aligned}
\trace_\fock q^{L_0} &\D(t_1) \cdots\D(t_n) \\
 &=
\trace_{\fock^\half\otimes\fock^\half} q^{L_0}(\D_1(t_1) +
\D_2(t_1))\cdots (\D_1(t_n) + \D_2(t_n)) \\
 &=
\sum_{\vec{i}\in\{1,2\}^n} \trace_{\fock^\half\otimes\fock^\half}
q^{L_0}\D_{i_1}(t_1)\D_{i_2}(t_2)\cdots\D_{i_n}(t_n).
\end{aligned}$$
which is equivalent to the first formula in the proposition.

Observe that $$\begin{aligned} \alpha' &= \sum_{k>0}\left(\varphi_{-k}\varphi_{k} + \varphi_{-k}'\varphi_k'\right) \\ &= \sum_{k>0}\left(\psi_{-k}^+\psi_{k}^- + \psi_{-k}^-\psi_{k}\right)\end{aligned}$$ from the isomorphism of Fock spaces.  Recalling that $$e_{11} = \sum_{k>0}\left(\psi_{-k}^+\psi_{k}^- - \psi_{-k}^-\psi_{k}\right)$$ it follows that $$\left(-1\right)^\alpha = \left(-1\right)^{e_{11}}.$$ Thus we have $$ \trace_\fock \left(-1\right)^{\alpha'} q^{L_0}\D(t_1)\cdots\D(t_n) = \trace_\fock \left(-1\right)^{e_{11}} q^{L_0}\D(t_1)\cdots\D(t_n);$$ the proposition follows by noting that on the right-hand side of (\ref{eq:subset}), there are exactly two terms equal to $\overline{\mathbb D}^\hf_{(0)}(q;{\bf t})$ which come from $I = \emptyset$ and $\{1,\dots,n\}$.  Note that a formula for $\trace_\fock \left(-1\right)^{e_{11}} q^{L_0}\D(t_1)\cdots\D(t_n)$ is given by (\ref{needed}) with $z = -1$.
\end{proof}

We now present our main theorem.

\begin{theorem} \label{th:lplushalftwist}
The function
$\overline{\mathbb D}^{l+\hf}_\la(q,{\bf t})$ is equal to
\begin{eqnarray*}
 && \overline{\mathbb D}^{\hf}_{(0)} (q;{\bf t}) \times \\
 &&\times
 \sum_{\sigma\in W(B_l)}\left(-1\right)^{\ell(\sigma)}
q^{\frac{\Vert \la+\rho-\sigma(\rho) \Vert^2}{2}}
\prod_{a=1}^l
 \Big(
\sum_{\vec{\eps}_a\in\{\pm 1\}^n}
 [\vec{\eps}_a] (\Pi {\bf t}^{\vec{\eps}_a})^{k_a}
 F_{bo}(q;{\bf t}^{{\vec{\eps}_a}})
 \Big)
\end{eqnarray*}
where $k_a = (\la+\rho-\sigma(\rho), \ep_a)$.
\end{theorem}

\begin{proof} From \cite[Theorem 4.1]{W1}, as $(O(2l+1),d_\infty)$-modules,
\begin{equation}\label{dualitystatement}\fock^{l+\half} \cong \bigoplus_{\la\in\Sigma(B)} V_\la (O(2l+1))
\otimes L(d_\infty, \Lambda(\la)).\end{equation}  Apply $\trace_{\fock^{l+\half}} \left(-1\right)^\alpha z_1^{e_{11}}\cdots z_l^{e_{ll}} q^{L_0}\D(t_1)\cdots\D(t_n)$ to both sides of this Howe-duality decomposition.  As $\alpha$ only acts on $\fock^{\half}$, we obtain $$\begin{aligned}\trace_{\fock^{\half}} \left(-1\right)^\alpha q^{L_0}\D(t_1)&\cdots\D(t_n)\cdot\prod_{i=1}^l \trace_{\fock} z_{i}^{e_{ii}}q^{L_0}\D(t_1)\cdots\D(t_n) \\& = \sum_{\la\in\Sigma(B)} \left(-1\right)^{\det+1} \mathrm{ch}^b_\la(z_1,\dots,z_l)\mathbb{D}_\la^{l+\half}(q,{\bf t})\end{aligned}$$ where $(-1)^{\det+1}$ is equal to 1 if $\la\in\PA^l$ and -1 otherwise.

For $\la\in \PA^l$, the character of the irreducible
$O(2l+1)$-module associated to $\la$ and $\la\otimes \det$ is the
same, and is given as follows (cf. \cite[p. 408]{FH})
\begin{eqnarray} \label{eq:charB}
\text{ch}^b_\la(z_1,\dots,z_l)
 = \frac{\left|z_j^{\la_i +l -i +\half} -z_j^{-(\la_i +l -i
+\half)} \right|}{\left|z_j^{l-i +\half} -z_j^{-(l -i
+\half)}\right|},
\end{eqnarray}
so we may rewrite the above as
$$\begin{aligned}\overline{\mathbb D}^{\hf}_{(0)}(q;{\bf t}) &\prod_{i=1}^l \trace_{\fock} z_{i}^{e_{ii}}q^{L_0}\D(t_1)\cdots\D(t_n) \\ &= \sum_{\la\in\Sigma(B)} \left(-1\right)^{\det+1} \frac{\left|z_j^{\la_i +l -i +\half} -z_j^{-(\la_i +l -i
+\half)} \right|}{\left|z_j^{l-i +\half} -z_j^{-(l -i +\half)}\right|} \mathbb{D}_\la^{l+\half}(q,{\bf t}).\end{aligned}$$
The Weyl denominator of type $B_l$ reads that
\begin{eqnarray} \label{eq:denomB}
 {\left |z_j^{l-i+\hf}+z_j^{-(l-i+\hf)}\right|} = \sum_{\sigma\in
W(B_l)} (-1)^{\ell(\sigma)} {\bf z}^{\sigma(\rho)}.
\end{eqnarray}
so by cross-multiplying terms, we may write
$$\begin{aligned}\sum_{\sigma\in W(B_l)} (-1)^{\ell(\sigma)} &{\bf z}^{\sigma(\rho)}\cdot\overline{\mathbb D}^{\hf}_{(0)}(q;{\bf t}) \prod_{i=1}^l \trace_{\fock} z_{i}^{e_{ii}}q^{L_0}\D(t_1)\cdots\D(t_n) \\ &= \sum_{\la\in\Sigma(B)} \left(-1\right)^{\det+1} \left|z_j^{\la_i +l -i +\half} -z_j^{-(\la_i +l -i
+\half)} \right| \mathbb{D}_\la^{l+\half}(q,{\bf t}).\end{aligned}$$
We may use (\ref{needed}) to expand each multiplicand on the left-hand side and compare coefficients of the dominant monomial $\mathbf{z}^{\la+\rho}$ on each side to finish the proof.
\end{proof}

\begin{remark} In the spirit of this paper, there are three more cases where one can consider the correlation functions on irreducible components of a direct summand.  At the positive level, the integral case of $\dinf$ remains.  This case is technically more difficult and different than the case we consider here; we do not know of an operator in $\dinf$ that is able to differentiate between the two modules of a direct summand, and this phenomenon only occurs for certain irreducible modules.  A much different strategy may be required.

Also, at the negative level, both the integral and half-integral cases for $\cinf$ are similar to $\dinf$ for the positive levels.  Given the already different nature of the negative level cases, again, a much different strategy may be required.  The author plans to consider these in the future.\end{remark}

\subsection{The $q$-dimension of a $\dinf$-module of level $l+\hf$}

For a $d_\infty$-module $M$, we denote by $$\dim_q M = \trace_M q^{L_0}$$ the $q$-dimension of the module $M$.  Set $$\mathsf{Q}(q)^+ =  \trace_{(L(d_\infty;\Lambda(\la)) \oplus
L(d_\infty;\Lambda({\la}\otimes\det))} q^{L_0}$$ and $$\mathsf{Q}(q)^- = \trace_{(L(d_\infty;\Lambda(\la)) \oplus
L(d_\infty;\Lambda({\la}\otimes\det))} \left(-1\right)^\alpha q^{L_0}.$$

The following proposition is a direct consequence of the above notation and the proof of Proposition \ref{maintheorem}.

\begin{prop}\label{qprop} For $\la\in\PA^l$, we have $$\dim_q L(d_\infty,\Lambda(\la)) = \frac{ \mathsf{Q}(q)^+ + \mathsf{Q}(q)^-}{2}$$ and
$$\dim_q L(d_\infty,\Lambda(\la\otimes\det)) = \frac{\mathsf{Q}(q)^+ - \mathsf{Q}(q)^-}{2}.$$
\end{prop}

Note that $\mathsf{Q}(q)^+ = \dim_q [L(d_\infty;\Lambda(\la)) \oplus
L(d_\infty;\Lambda({\la}\otimes\det))]$ which is equal to the following equivalent formulas (cf. \cite{TW}):

\begin{align*}
 \mathsf{Q}(q)^+ & =
\frac{(- q^{-\hf};q)_\infty}{(q;q)_\infty^l} \cdot
\sum_{\sigma\in W(B_l)}(-1)^{\ell(\sigma)} q^{\frac{\Vert
\la+\rho-\sigma(\rho) \Vert^2}{2}} \\
 &=
\frac{(- q^{-\hf};q)_\infty}{(q;q)_\infty^l} \cdot
q^{\frac{\Vert\la\Vert^2}{2}}
  \prod_{1\leq i\leq l} \left(1-q^{\la_i+l-i+1/2}\right) \times\\
 & \quad \times \prod_{1\leq i<j\leq l}
\left(1-q^{\la_i-\la_j+j-i}\right)
\left(1-q^{\la_i+\la_j+2l-i-j+1}\right).
\end{align*}

It remains to compute $\mathsf{Q}(q)^-$.

\begin{prop}
We have
\begin{align*}
 \mathsf{Q}(q)^- & =
\frac{(q^{-\hf};q)_\infty}{(q;q)_\infty^l} \cdot
\sum_{\sigma\in W(B_l)}(-1)^{\ell(\sigma)} q^{\frac{\Vert
\la+\rho-\sigma(\rho) \Vert^2}{2}} \\
 &=
\frac{(q^{-\hf};q)_\infty}{(q;q)_\infty^l} \cdot
q^{\frac{\Vert\la\Vert^2}{2}}
  \prod_{1\leq i\leq l} \left(1-q^{\la_i+l-i+1/2}\right) \times\\
 & \quad \times \prod_{1\leq i<j\leq l}
\left(1-q^{\la_i-\la_j+j-i}\right)
\left(1-q^{\la_i+\la_j+2l-i-j+1}\right).
\end{align*}
\end{prop}

\begin{proof} In the proof of Theorem \ref{th:lplushalftwist}, we instead apply $$\trace_{\fock^{l+\hf}} \left(-1\right)^\alpha z_1^{e_{11}}\cdots z_l^{e_{ll}} q^{L_0}$$ to both sides of the duality in (\ref{dualitystatement}).  The same strategy applies, with the substitutions $$\trace_{\fock^\hf} \left(-1\right)^\alpha q^{L_0} = (q^{-\hf};q)_\infty$$ and $$\trace_{\fock} z_i^{e_{ii}} q^{L_0} = \dim_q \fock^(0) \sum_{k\in\Z} z_i^k q^{k^2/2} = (q;q)^{-1}_\infty \sum_{k\in\Z} z_i^k q^{\frac{k^2}{2}}.$$ The equivalence of the two statements follows from above.\end{proof}

We note that the $q$-dimension formula was also obtained in an alternate form using a very different strategy in \cite{KWY}.

\end{document}